# NEW APPROACHES TO BAYESIAN CONSISTENCY[1]

### By Stephen Walker

#### *University of Bath*

We use martingales to study Bayesian consistency. We derive sufficient conditions for both Hellinger and Kullback–Leibler consistency, which do not rely on the use of a sieve. Alternative sufficient conditions for Hellinger consistency are also found and demonstrated on examples.

**1. Introduction.** Let $X_1, X_2, \ldots$, taking values in $(\mathscr{X}, \mathscr{B})$, be independent and identically distributed random variables from some fixed but unknown (the true) density function $f_0$, with corresponding distribution function $F_0$. Let $F_0^n$ be the corresponding $n$-fold product measure on $(\mathscr{X}^n, \mathscr{B}^n)$ and let $F_0^\infty$ denote the infinite product measure.

With $f_0$ being unknown, the Bayesian constructs a prior distribution $\Pi$ on $\Omega$, the space of density functions on $(\mathscr{X}, \mathscr{B})$. This prior combines with the data to define the posterior distribution $\Pi^n$, assigning mass

$$(1) \qquad \Pi^n(A) = \frac{\int_A R_n(f) \Pi(df)}{\int R_n(f) \Pi(df)}$$

to the set of densities $A$, where

$$R_n(f) = \prod_{i=1}^n f(X_i)/f_0(X_i).$$

The predictive density is given by

$$f_n(x) = \int f(x) \Pi^n(df).$$

Here, and throughout, we assume that all relevant $f$ are, in fact, densities with respect to the Lebesgue measure.

---

Received April 2003; revised January 2004.

[1]Supported by an EPSRC Advanced Research Fellowship.

---

This paper is concerned with Hellinger and Kullback–Leibler consistency. For example, for Hellinger consistency, the required result we are aiming for is

$$\Pi^n(\{f : H(f, f_0) > \varepsilon\}) \to 0 \qquad \text{a.s. } [F_0^\infty]$$

for all $\varepsilon > 0$, where $H(f, f_0)$ is the Hellinger distance between $f$ and $f_0$, given by

$$H(f, f_0) = \left\{ \int (\sqrt{f} - \sqrt{f_0})^2 \right\}^{1/2}.$$

Previous studies of Hellinger consistency [see, e.g., Barron, Schervish and Wasserman (1999) and Ghosal, Ghosh and Ramamoorthi (1999)] deal with the numerator and denominator in the expression for $\Pi^n(A)$ separately. Briefly, if $\Pi$ puts positive mass on all Kullback–Leibler neighborhoods of $f_0$ (which will be referred to as the Kullback–Leibler property for $\Pi$), then the denominator is eventually bounded below by $\exp(-nc)$ for all $c > 0$. Setting $A = \{f : H(f, f_0) > \varepsilon\}$, for some $\varepsilon > 0$, with constraints on the prior, ensuring sufficiently low mass on densities which track the data too closely, the numerator can be eventually bounded above by $\exp(-n\delta)$, for some $\delta > 0$. Consequently, with the appropriate conditions in place, $\Pi^n(A) \to 0$ a.s., with exponential rate, for all $\varepsilon > 0$.

To be more explicit, the basic ideas of current approaches are based on the introduction of an increasing sequence of sets $\mathscr{G}_n$, a sieve, and to consider

$$\Pi^n(A) = \Pi^n(\mathscr{G}_n \cap A) + \Pi^n(\mathscr{G}_n^c \cap A).$$

Putting sufficiently low mass on densities which track the data too closely, that is, the densities in $\mathscr{G}_n^c$, involves ensuring that $\Pi(\mathscr{G}_n^c) < \exp(-n\xi)$ for all large $n$ and for some $\xi > 0$. This results in $\Pi^n(\mathscr{G}_n^c) < \exp(-n\xi^*)$ a.s. for all large $n$ for some $\xi^* > 0$. The aim then is to find $\mathscr{G}_n$ such that

$$\int_{A \cap \mathscr{G}_n} R_n(f) \Pi(df) < \exp(-n\delta) \qquad \text{a.s.}$$

for all large $n$ for some $\delta > 0$. Approaches differ in the precise form of $\mathscr{G}_n$ which guarantees the above. For example, Ghosal, Ghosh and Ramamoorthi (1999) have $J(\eta, \mathscr{G}_n) < n\beta_\eta$ for all large $n$, for some $\beta_\eta > 0$ for all $\eta > 0$, where $J$ is the $L_1$ metric entropy.

We also deal with the numerator and denominator separately but study the numerator via different techniques which include the use of martingales. We do not use sieves. To fix the notation, define

$$f_{nA}(x) = \int f(x) \Pi_A^n(df)$$

to be the predictive density with posterior distribution restricted, and normalized, to the set $A$, let $h(f, f_0) = 1 - \int \sqrt{f f_0} = H^2(f, f_0)/2$ be a slight variation on the Hellinger distance $H$, and note that $h(f, f_0) \leq 1$. Also define $I_n = \int R_n(f)\Pi(df)$ and $D(f, f_0) = \int f_0 \log(f_0/f)$ to be the Kullback–Leibler divergence between $f$ and $f_0$. The Kullback–Leibler property is given by

$$\Pi(\{f : D(f, f_0) < \varepsilon\}) > 0$$

for all $\varepsilon > 0$. Since $f_0$ is unknown, the condition is

$$\Pi(\{f : D(f, g) < \varepsilon\}) > 0$$

for all densities $g$ and all $\varepsilon > 0$. This is possible to achieve using nonparametric priors. See Barron, Schervish and Wasserman (1999) and Ghosal, Ghosh and Ramamoorthi (1999).

The layout of the paper is as follows. In Section 2 we present preliminary results based on certain martingales. Section 3 unifies approaches to posterior consistency via the use of these martingales and Section 4 deals with the special case of consistency for predictive densities. Section 5 presents a specific result for Hellinger consistency which does not use martingales and examples are presented in Section 6. Section 7 contains a discussion and highlights areas for future research.

**2. Preliminaries.** Here we will discuss fundamental concepts and ideas on which the paper is based. Our concern is with the numerator $L_n = \int_A R_n(f)\Pi(df)$, where $A$ is a set of densities, of (1). We have already established that the Kullback–Leibler property will always deal appropriately with the denominator. The following identity is the key:

$$(2) \qquad L_{n+1}/L_n = f_{nA}(X_{n+1})/f_0(X_{n+1}), \qquad n = 0, 1, \ldots,$$

and it is easy to check that this holds. From here we can go in one of two directions. The first option is based on martingales and takes $A = \{f : d(f, f_0) > \varepsilon\}$, where $d$ metricizes weak convergence, and is the Hellinger or the Kullback–Leibler distance. The second option, in the case of the Hellinger distance, is to split $\{f : H(f, f_0) > \varepsilon\}$ into a countable number of disjoint sets $\{A_j\}$ based on Hellinger balls, $A_j = \{f : H(f, f_j) < \delta\}$ for some suitable set of densities $\{f_j\}$ and some $\delta > 0$. This is possible due to the separability of $\Omega$ with respect to the Hellinger metric.

The two approaches share similarities, both use (2), but are otherwise different. The first covers a range of types of consistency, whereas the second seems suited only to Hellinger consistency. To set the scene for the first option we consider measurable functions $T_d$, linked to a distance measure $d$, such that

$$\mathrm{E}\{T_d(L_{n+1}/L_n)|\mathscr{F}_n\} = -d(f_{nA}, f_0),$$

where $\mathscr{F}_n = \sigma(X_1, \ldots, X_n)$. If $T_d(y) = \sqrt{y} - 1$, then $d(f, f_0) = h(f, f_0)$ and if $T_d(y) = \log y$, then $d(f, f_0) = D(f, f_0)$. Other cases arise; for example, if $T_d(y) = 1 - 1/y$, then $d(f, f_0) = \int f_0^2/f - 1$, which is the $\chi$-squared distance.

Now consider the martingale $(M_N, \mathscr{F}_N)$ given by

$$M_N = \sum_{n=1}^{N} \{T_d(L_n/L_{n-1}) + d(f_{n-1A}, f_0)\}.$$

A well-known result for such martingales [see Loève (1963)] is that if

(3)                          $$\sum_n n^{-2} \operatorname{Var}\{T_d(L_n/L_{n-1})\} < \infty,$$

then $M_N/N \to 0$ a.s. Consequently, if

$$\liminf_n d(f_{nA}, f_0) > 0 \qquad \text{a.s.},$$

then

$$\limsup_N \frac{1}{N} \sum_{n=1}^{N} T_d(L_n/L_{n-1}) < 0 \qquad \text{a.s.}$$

For both cases of $T_d(y) = \sqrt{y} - 1$ and $T_d(y) = \log y$, the above implies that there exists a $\delta > 0$ such that $L_N < \exp(-N\delta)$ a.s. for some $\delta > 0$ for all large $N$. This result can be achieved for $T_d(y) = \sqrt{y} - 1$ by making use of the fact that an arithmetic mean is greater than or equal to a geometric mean, and it is clearly true for $T_d(y) = \log y$. It is worth writing this down formally.

LEMMA 1.    *Let $L_n = \int_A R_n(f)\Pi(df)$ and $T_d(y) = \sqrt{y} - 1$ or $T_d(y) = \log y$. If* (3) *holds and*

$$\liminf_n d(f_{nA}, f_0) > 0 \qquad a.s.,$$

*then $L_N < \exp(-N\delta)$ a.s. for some $\delta > 0$ for all large $N$.*

This result, namely, $L_N < \exp(-N\delta)$ a.s. for some $\delta > 0$ for all large $N$, combined with the Kullback–Leibler property for $\Pi$, leads to $\Pi^n(A) \to 0$ a.s. This follows since the Kullback–Leibler property implies $I_N > \exp(-Nc)$ a.s. for all large $N$, for any $c > 0$. Hence, we can choose $c < \delta$.

**3. Posterior consistency.**    In this section we unify posterior consistency based on Lemma 1. Here we will drop the subscript $d$ from $T$.

3.1. *Weak consistency.*    Here we have $A = \{f : d_W(f, f_0) > \varepsilon\}$, where $d_W$ metricizes weak convergence of probability distributions, that is, $d_W(f_n, f_0) \to 0$ if and only if $\int g(x) f_n(x)\, dx \to \int g(x) f_0(x)\, dx$ for all continuous and bounded $g$. Now $H(f_{nA}, f_0) > \varepsilon^*$ for all large $n$ a.s. for some $\varepsilon^* > 0$ since eventually $f_{nA}$ does not lie in a weak neighborhood of $f_0$ and so neither does it lie in a Hellinger neighborhood of $f_0$. Hence, taking $T(y) = \sqrt{y} - 1$, we have

$$(4) \qquad \sum_n n^{-2} \operatorname{Var}\{T(L_n/L_{n-1})\} < \infty$$

automatically as $\mathrm{E}(L_n/L_{n-1}) = 1$. Hence, both conditions of Lemma 1 are satisfied and so the Kullback–Leibler property is sufficient for weak consistency. This is, of course, known; see Schwartz (1965).

3.2. *Hellinger consistency.* Here we retain $T(y) = \sqrt{y} - 1$ and consider $A = \{f : H(f, f_0) > \varepsilon\}$. While, as in Section 3.1, we remain with (4) being true, we do not automatically have $\liminf_n H(f_{nA}, f_0) > 0$ a.s. Hence, we have only one of the conditions of Lemma 1 being satisfied automatically.

THEOREM 1. *If $\Pi$ has the Kullback–Leibler property, then*

$$\Pi^n(A) \to 0 \qquad a.s.$$

*for all sets $A$ for which $\liminf_n H(f_{nA}, f_0) > 0$ a.s.*

This extends Walker (2003) who showed that if $\Pi$ has the Kullback–Leibler property and $H(f_{nA}, f_0) > \gamma$ for all $n$ a.s. for some $\gamma > 0$, then $\Pi^n(A) \to 0$ a.s. This result was then used to obtain the Hellinger consistency result of Ghosal, Ghosh and Ramamoorthi (1999).

3.3. *Kullback–Leibler consistency.* In view of the importance of the Kullback–Leibler property to Bayesian consistency, it would make sense to find additional sufficient conditions for posteriors to accumulate in all Kullback–Leibler neighborhoods of $f_0$. There are also practical reasons. A Bayesian approach to parametric prediction advocated by Walker and Gutiérrez-Peña (1999) entails minimizing $D(f_n, f_\lambda)$. Here $f_\lambda$ is a parametric family of densities and $f_n$ is a nonparametric predictive density. For large sample suitability of this procedure it is important that $D(f_n, f_0) \to 0$ a.s. Further motivation for Kullback–Leibler consistency is given in Barron (1988) who cites universal data compression and stock market portfolio selection as applications where this type of consistency is important.

For the martingale $M_N$ we now take $T(y) = \log y$ and consider $A = \{f : D(f, f_0) > \varepsilon\}$. In this case neither of the conditions of Lemma 1 holds automatically.

THEOREM 2. *If $\Pi$ has the Kullback–Leibler property and*

$$(5) \qquad \sum_n n^{-2} \operatorname{Var}\{\log(L_n/L_{n-1})\} < \infty,$$

*then $\Pi^n(A) \to 0$ a.s. for all sets $A$ for which $\liminf_n D(f_{nA}, f_0) > 0$ a.s.*

To examine Theorem 2 further, we write $L_n = \Pi^n(A)I_n$, giving

$$M_N = \log I_N + \log\{\Pi^N(A)/\Pi(A)\} + \sum_{n=1}^{N} D(f_{n-1A}, f_0).$$

If $\Pi$ has the Kullback–Leibler property, then $N^{-1}\log I_N \to 0$ a.s. This follows since $I_n > \exp(-nc)$ a.s. for all large $n$ for any $c > 0$ and, because $E(I_n) = 1$, we have $I_n < \exp(nc)$ a.s. for all large $n$ for any $c > 0$. So, $M_N/N \to 0$ a.s. and $\Pi^N(A) < \exp(-Nc)$ a.s. for some $c > 0$ for all large $N$ together imply that

$$(6) \qquad \liminf_N \frac{1}{N}\sum_{n=1}^{N} D(f_{n-1A}, f_0) > c \qquad \text{a.s.}$$

Hence, Theorem 2 could be written as follows:

THEOREM 2*. *If $\Pi$ has the Kullback–Leibler property and* (5) *holds, then* $\Pi^n(A) < \exp(-nc)$ *a.s. for all large $n$ if and only if* (6) *holds.*

If $A = \{f : D(f, f_0) > \varepsilon\}$, then one anticipates that $\liminf_n D(f_{nA}, f_0) \geq \varepsilon$ a.s. However, it is difficult to establish when $\liminf_N N^{-1}\sum_{n=1}^{N} D(f_{n-1A}, f_0) > 0$ a.s., yet, when $\Pi$ has the Kullback–Leibler property and (5) holds, which is not a particularly demanding condition, it does become a necessary condition for Kullback–Leibler consistency with exponential rate. It should also be pointed out that Theorem 2* equally applies to Hellinger consistency when $A = \{f : H(f, f_0) > \varepsilon\}$ and the necessary condition also applies.

**4. Predictive consistency.** Here we take $A = \Omega$ so that $f_{nA} \equiv f_n$, the predictive density. Also, $L_n \equiv I_n$, the denominator of (1). Hence,

$$M_N = \sum_{n=1}^{N}\{T(I_n/I_{n-1}) + d(f_{n-1}, f_0)\}.$$

LEMMA 2. *If $\Pi$ has the Kullback–Leibler property and $T(y) = \sqrt{y} - 1$ or $T(y) = \log y$, then*

$$\frac{1}{N}\sum_{n=1}^{N} T(I_n/I_{n-1}) \to 0 \qquad a.s.$$

PROOF. This is obvious with $T(y) = \log y$ since $N^{-1}\log I_N \to 0$ a.s. as $N \to \infty$ when $\Pi$ has the Kullback–Leibler property. When $T(y) = \sqrt{y} - 1$, we know that $M_N/N \to 0$ a.s. and so

$$\frac{1}{N}\sum_{n=1}^{N} T(I_n/I_{n-1}) + \frac{1}{N}\sum_{n=1}^{N} h(f_{n-1}, f_0) \to 0 \qquad \text{a.s.}$$

Now

$$\frac{1}{N} \sum_{n=1}^{N} T(I_n/I_{n-1}) \geq I_N^{1/(2N)} - 1 \to 0 \qquad \text{a.s.}$$

and so

$$\frac{1}{N} \sum_{n=1}^{N} T(I_n/I_{n-1}) \to 0 \qquad \text{a.s.}$$

as $h \geq 0$, completing the proof. $\square$

If we take $T(y) = \sqrt{y} - 1$, then $M_N/N \to 0$ a.s. and, from Lemma 2, we have

$$\frac{1}{N} \sum_{n=1}^{N} H(f_{n-1}, f_0) \to 0 \qquad \text{a.s.}$$

This result is found in Walker (2003). The following theorem applies by considering $T(y) = \log y$.

THEOREM 3. *If $\Pi$ has the Kullback–Leibler property and*

(7)
$$\sum_n n^{-2} \operatorname{Var}\{\log(I_n/I_{n-1})\} < \infty,$$

*then*

$$\frac{1}{N} \sum_{n=1}^{N} D(f_{n-1}, f_0) \to 0 \qquad a.s.$$

It is straightforward to demonstrate that (7) holds when

$$\sup_n \left\{ \mathrm{E}_{X^n} \int f_0^2/f_n \right\} < \infty.$$

Here $\mathrm{E}_{X^n}$ is the expectation with respect to $X^n = (X_1, \ldots, X_n)$ taken independently from $f_0$. See Section 6.4 for an example illustrating a non-parametric prior for which $\sup_n \{\mathrm{E}_{X^n} \int f_0^2/f_n\} < \infty$.

**5. Hellinger consistency.** To introduce the ideas here, consider the discrete prior which puts mass $\Pi_k$ on the density function $f_k$, for $k = 1, 2, \ldots$. In this case the posterior mass assigned to $f_k$ is given by

$$\Pi_k^n = \frac{R_n(f_k)\Pi_k}{\sum_j R_n(f_j)\Pi_j}.$$

If we assume the prior has the Kullback–Leibler property, then the additional condition for Hellinger consistency turns out to be

$$\sum_k \sqrt{\Pi_k} < \infty.$$

REMARK 1.   The result provides information concerning the counterexample appearing in Barron, Schervish and Wasserman (1999), which shows that the Kullback–Leibler property for $\Pi$ is not sufficient for Hellinger consistency. The prior in this case puts positive mass on single densities and, for each integer $N$, has sets of these densities $\mathscr{P}_N$ for which $\Pi(\mathscr{P}_N) > \eta/N^2$ for some $\eta > 0$. Clearly, then

$$\sum_N \sqrt{\Pi(\mathscr{P}_N)} = \infty.$$

Now $\Omega$ is separable; that is, we can cover $\Omega$ with a countable set of Hellinger balls of radius $\delta$ for any $\delta > 0$. Therefore,

$$A = \{f : h(f, f_0) > \varepsilon\}$$

can be covered by the countable union of disjoint sets $A_j$, where $A_j \subseteq A_j^* = \{f : h(f, f_j) < \delta\}$, and $\{f_j\}$ is a set of densities such that $h(f_j, f_0) > \varepsilon$. We can take $\delta < \varepsilon$ so that $h(f_{nA_j}, f_0) > \gamma > 0$, where $\gamma = \varepsilon - \delta$. This follows since

$$h(f_{nA_j}, f_0) \geq h(f_j, f_0) - h(f_{nA_j}, f_j)$$

and $h(f_{nA_j}, f_j) \leq \delta$.

THEOREM 4.   *If $\Pi$ has the Kullback–Leibler property and*

$$\sum_j \sqrt{\Pi(A_j)} < \infty,$$

*then*

$$\Pi^n(A) \to 0 \qquad a.s.$$

PROOF.   Now

$$\Pi^n(A) = \sum_j \Pi^n(A_j) \leq \sum_j \sqrt{\Pi^n(A_j)}$$

$$= \sum_j \sqrt{\int_{A_j} R_n(f)\Pi(df)/I_n}.$$

If

$$\Lambda_{nj} = \sqrt{\int_{A_j} R_n(f)\Pi(df)},$$

then

$$\Lambda_{n+1j} = \Lambda_{nj}\sqrt{f_{nA_j}(X_{n+1})/f_0(X_{n+1})};$$

see Section 2, equation (2), which includes the case when $n = 0$ and $\Lambda_{0j} = \sqrt{\Pi(A_j)}$. So,

$$\mathrm{E}(\Lambda_{n+1j}|\mathscr{F}_n) = \Lambda_{nj}\{1 - h(f_{nA_j}, f_0)\} < (1 - \gamma)\Lambda_{nj}$$

and, hence, $\mathrm{E}(\Lambda_{nj}) < (1 - \gamma)^n\sqrt{\Pi(A_j)}$. Therefore,

$$\mathrm{pr}\left\{\sum_j \Lambda_{nj} > \exp(-nd)\right\} < \exp(nd)(1 - \gamma)^n \sum_j \sqrt{\Pi(A_j)}$$

and so if

$$\sum_j \sqrt{\Pi(A_j)} < \infty,$$

then

$$\sum_j \Lambda_{nj} < \exp(-nd) \qquad \text{a.s.}$$

for all large $n$, for any $d < -\log(1 - \gamma)$. The Kullback–Leibler property for $\Pi$ ensures that $I_n > \exp(-nc)$ a.s. for all large $n$, for any $c > 0$. This completes the proof. $\square$

Clearly, if the prior $\Pi$ puts mass $\Pi_k$ on the density $f_k$ for $k = 1, 2, \ldots,$ then the required condition is simply

$$\sum_k \sqrt{\Pi_k} < \infty,$$

which is straightforward to arrange in practice.

The result of Theorem 4 can be applied to specific priors with good results. See next in Section 6. However, it does somewhat lack interpretation as can be seen by the need to go from $\sum_j \Pi^n(A_j)$ to $\sum_j \sqrt{\Pi^n(A_j)}$. On the other hand, the appearance of square roots should not be a great surprise when dealing with the Hellinger distance.

**6. Illustrations.** Here we consider some examples (6.1 to 6.3) illustrating Theorem 4. We have $\Omega$ being covered by $\{A_1, A_2, \ldots\}$, which are mutually disjoint Hellinger balls of radius $\delta$. The aim then is to show that

$$\sum_k \sqrt{\Pi(A_k)} < \infty$$

and that this holds for all $\delta > 0$. Also, Section 6.4 will illustrate Theorem 3.

6.1. *Infinite-dimensional exponential families.*   Here we consider the case when $f$ is constructed from an infinite sequence of random variables, $\theta_1, \theta_2 \dots$. The prior on the $\{\theta_j\}$ makes them independent and we assume that $\theta_j \sim \mathrm{N}(0, \sigma_j^2)$. A $\delta$-covering of $\Omega$ will be the union of sets of the type

$$\{\theta : n_j \delta_j < \theta_j < (n_j + 1)\delta_j, j = 1, 2, \dots\}$$

for a sequence $\delta_j = \delta \gamma_j$, where the $\{\gamma_j\}$ do not depend on $\delta$. Here the $n_j$ are integers and can be between $-\infty$ and $+\infty$. It is convenient to define

$$A_{jn} = (n\delta_j, (n+1)\delta_j).$$

We are then interested in the finiteness of

$$\sum_{n_1 = -\infty}^{\infty} \cdots \sum_{n_M = -\infty}^{\infty} \prod_{j=1}^{M} \sqrt{\mathrm{pr}(\theta_j \in A_{jn_j})}$$

as $M \to \infty$, which, because of symmetry, holds if

$$\prod_{j=1}^{\infty} \sum_{n=0}^{\infty} \sqrt{\mathrm{pr}(\theta_j \in A_{jn})} < \infty.$$

Dropping the subscript $j$ temporarily, we have

$$\sum_{n=0}^{\infty} \sqrt{\mathrm{pr}(\theta \in A_n)} \leq 1 + \sum_{n=1}^{\infty} \sqrt{\mathrm{pr}(\theta \in A_n)}$$

$$\leq 1 + (2\pi)^{-1/4} (\delta/\sigma)^{1/2} \sum_{n=1}^{\infty} \exp\{-\delta^2 n^2/(4\sigma^2)\}$$

$$\leq 1 + (2\pi)^{-1/4} (\delta/\sigma)^{1/2} [\exp\{\delta^2/(4\sigma^2)\} - 1]^{-1}$$

$$\leq 1 + 4^m m! (2\pi)^{-1/4} (\sigma/\delta)^{2m-1/2}$$

for any $m = 1, 2, \dots$. The last inequality follows from

$$\frac{\xi^{1/4}}{e^{\xi/4} - 1} \leq 4^m m! \xi^{1/4-m}, \qquad m = 1, 2, \dots,$$

for all $\xi > 0$. The required condition on the $\{\sigma_j\}$ is then that

$$\prod_{j=1}^{\infty} \{1 + \psi(\sigma_j/\gamma_j)^{2m-1/2}\} < \infty$$

for all $\psi > 0$. This is achieved if

$$\sum_j (\sigma_j/\gamma_j)^{2m-1/2} < \infty.$$

To make this example specific, consider the infinite-dimensional exponential family on $[0, 1]$ for which

$$f(x) = \exp\left\{\sum_{j=0}^{\infty} \theta_j \phi_j(x) - c(\Theta)\right\},$$

where the $\{\phi_j\}$ are an orthonormal basis on $[0, 1]$ and $c(\Theta)$ ensures that $f$ integrates to 1. Such an orthonormal basis is given by

$$\phi_0(x) = 1 \quad \text{and} \quad \phi_j(x) = \sqrt{2}\cos(j\pi x) \qquad \text{for } j \geq 1.$$

To ensure that $f$ is a density with probability 1, it is sufficient that $\sum_j \sigma_j < \infty$. Then, according to Barron, Schervish and Wasserman (1999), an additional condition sufficient for Hellinger consistency is that $\sum_j j\sigma_j < \infty$. So it is possible to have $\sigma_j \propto j^{-2-r}$ for any $r > 0$.

For our condition, we require

$$\sum_j (\sigma_j/\omega_j)^{2m-1/2} < \infty$$

for some sequence $\{\omega_j\}$ satisfying $\sum_j \omega_j < \infty$. This follows because we can take $\delta_j = \delta^* \omega_j / \sum_j \omega_j$ so that if $|\theta_{1j} - \theta_{2j}| < \delta_j$, then

$$\sup_{0 \leq x \leq 1} \left|\sum_j \theta_{1j}\phi_j(x) - \sum_j \theta_{2j}\phi_j(x)\right| < \delta^*\sqrt{2},$$

which implies that $h(f_1, f_2) < \delta = 1 - \exp(-\delta^*\sqrt{2})$. If we put $\omega_j \propto j^{-1-r}$ for any $r > 0$, then

$$\sum_j (\sigma_j j^{1+r})^{2m-1/2} < \infty$$

is sufficient. Therefore, we can actually have $\sigma_j \propto j^{-1-q}$ for any $q > 0$, by choosing $m$ large enough. This then is seen to be an improvement on the condition provided by Barron, Schervish and Wasserman (1999).

See also Walker and Hjort (2001) who have essentially the same result of $\sigma^j \propto j^{-1}$ as being sufficient for Hellinger consistency when combined with the Kullback–Leibler property for $\Pi$.

6.2. *Pólya trees.* Here we consider a Pólya tree prior on $[0, 1]$ with partition the dyadic intervals. Denote the sets at level $k$ by $B_{kj}$ for $j = 1, \ldots, 2^k$. Over these sets we have independent variables $\theta_{kj} \sim \text{be}(a_k, a_k)$ for odd $j$ and $\theta_{kj+1} = 1 - \theta_{kj}$, again for odd $j$. If $\sum_k a_k^{-1} < \infty$, then a random density function, with respect to the Lebesgue measure on $[0, 1]$, from the prior can be obtained via

$$f(x) = \lim_{k \to \infty} 2^k \prod_{j=1}^{k} \theta_{kj(x)},$$

where $kj(x)$ describes the interval at level $k$ in which $x$ lies; that is, $x \in B_{kj(x)}$. If $\sum_k a_k^{-1/2} < \infty$, then the prior puts positive mass on all Kullback–Leibler neighborhoods of densities $g$ for which $\int g \log g < \infty$. See, for example, Barron, Schervish and Wasserman (1999). However, according to Barron, Schervish and Wasserman (1999), the best sufficient condition for Hellinger consistency is $a_k = 8^k$ which is not a "nice" set-up and would impede those statisticians looking to incorporate relevant information. We improve on this using the sums of square-roots of prior probabilities.

First, we find the covering sets of $\Omega$ with Hellinger balls of radius $\delta$. If

$$\exp(-\delta_k) < \theta_{1kj}/\theta_{2kj} < \exp(\delta_k),$$

for all $j$, which is equivalent to

$$\exp(-\delta_k) < \theta_{1kj}/\theta_{2kj} < \exp(\delta_k)$$

and

$$\exp(-\delta_k) < (1 - \theta_{1kj})/(1 - \theta_{2kj}) < \exp(\delta_k)$$

for all odd $j$, and for all $k$ and $\sum_k \delta_k = \delta^*$, then it is easy to show that

$$\int \sqrt{f_1 f_2} > \exp(-\tfrac{1}{2}\delta^*)$$

and so $h(f_1, f_2) < \delta = 1 - \exp(-\delta^*/2)$. Here, for example,

$$f_1(x) = \lim_{k \to \infty} 2^k \prod_{l=1}^{k} \theta_{1lj(x)}.$$

Therefore, letting $\theta_k$ denote a generic random variable at level $k$, we split $[0,1]$, the range of $\theta_k$, into the sets $A_{k0} = (\tfrac{1}{2} - b_k, \tfrac{1}{2} + b_k)$, where

$$b_k = \tfrac{1}{2}\{\exp(\delta_k) - 1\}/\{\exp(\delta_k) + 1\}$$

and

$$A_{kl}^- = (c_k \exp\{-l\delta_k\}, c_k \exp\{-(l-1)\delta_k\})$$

and

$$A_{kl}^+ = (1 - c_k \exp\{-(l-1)\delta_k\}, 1 - c_k \exp\{-l\delta_k\})$$

for $l = 1, 2, \ldots$, where $c_k = \tfrac{1}{2} - b_k$. Again, due to symmetry, that is, $\mathrm{pr}(\theta_k \in A_{kl}^+) = \mathrm{pr}(\theta_k \in A_{kl}^-)$, we are interested in the finiteness, as $M \to \infty$, of

$$\sum_{n_{11}=0}^{\infty} \cdots \sum_{n_{2^M 2^{M-1}}=0}^{\infty} \prod_{k=1}^{M} \prod_{j=1,3,\ldots,2^k-1} \sqrt{\mathrm{pr}(\theta_{kj} \in A_{kn_{kj}}^-)},$$

with the convention that $A_{k0}^- = A_{k0}$, which is equivalent to the finiteness of

$$\prod_{k=1}^{\infty} \left\{ \sqrt{\mathrm{pr}(\theta_k \in A_{k0})} + \sum_{l=1}^{\infty} \sqrt{\mathrm{pr}(\theta_k \in A_{kl}^-)} \right\}^{2^{k-1}}.$$

The difference between this and the corresponding expression for the infinite-dimensional families is the power $2^{k-1}$, which is present due to the $2^{k-1}$ independent variables at level $k$. Now

$$\mathrm{pr}(\theta_k \in A_{kl}^-) = \frac{\Gamma(2a_k)}{\Gamma(a_k)^2} \int_{A_{kl}^-} x^{a_k-1}(1-x)^{a_k-1}\,dx$$

$$\leq 2^{2a_k-1} \sqrt{a_k/\pi}\, c_k \exp(-l\delta_k)\{\exp(\delta_k)-1\}\{\xi_{kl}(1-\xi_{kl})\}^{a_k-1},$$

where $\xi_{kl} = c_k \exp\{-(l-1)\delta_k\}$. Here we have used the inequality

$$\Gamma(2a)/\Gamma(a)^2 \leq 2^{2a-1}\sqrt{a/\pi},$$

see Barron, Schervish and Wasserman (1999), and that if $x < \xi < \frac{1}{2}$, then $x(1-x) < \xi(1-\xi)$. Now $\xi_{kl}(1-\xi_{kl}) \leq 1/4 - b_k^2$ for all $l$ and so

$$\sqrt{\mathrm{pr}(\theta_k \in A_{kl}^-)} \leq \psi 2^{a_k-1} a_k^{1/4} \exp(-\tfrac{1}{2}l\delta_k)\sqrt{\exp(\delta_k)-1}(\tfrac{1}{4}-b_k^2)^{a_k/2-1/2}$$

for some fixed $\psi > 0$. Here we need only consider $k$ large enough for which $a_k > 1$. Therefore,

$$\sqrt{\mathrm{pr}(\theta_k \in A_{k0})} + \sum_{l=1}^{\infty}\sqrt{\mathrm{pr}(\theta_k \in A_{kl}^-)} \leq 1 + \frac{\psi\sqrt{\exp(\delta_k)-1}}{\exp(\delta_k/2)-1} a_k^{1/4}(1-4b_k^2)^{a_k/2-1/2}.$$

We are then looking for

$$\sum_k 2^{k-1} a_k^{1/4} \frac{\sqrt{\exp(\delta_k)-1}}{\exp(\delta_k/2)-1} \exp(-2a_k b_k^2) < \infty.$$

Now we can take $\delta_k \propto k^{-1-r}$ for any $r > 0$ and for large $k$,

$$\frac{\sqrt{\exp(\delta_k)-1}}{\exp(\delta_k/2)-1} \sim \delta_k^{-1/2} = k^{1/2+r/2}.$$

Also, for large $k$, $b_k \sim \delta_k$ and so $a_k \propto k^{3+q}$ for any $q > 0$ is sufficient.

6.3. *Mixture of priors.* A popular type of nonparametric prior consists of a mixture of priors,

$$\Pi = \sum_N p_N \Pi_N,$$

where $\sum_N p_N = 1$ and the $\{p_N\}$ are fixed. Here $\Pi_N$ is supported by densities in $C_N \subseteq \Omega$, so that $\Pi_N(C_N) = 1$. Typically, $C_N$ is totally bounded, that is, $C_N$ can be covered by a finite number of Hellinger balls of size $\delta$, for any $\delta > 0$. The number of such balls will be denoted by $I_N(\delta)$. More often than not, $C_N \subseteq C_{N+1}$ and $C_N$ converges to $\Omega$. See Petrone and Wasserman (2002), for example, who consider random densities generated via Bernstein polynomials.

Following the above specifications, if $\bigcup_k A_k$ covers $\Omega$, then we assume that $C_N$ is covered by $\{A_1, \ldots, A_{I_N}\}$ and, therefore, $\Pi_N(A_k) = 0$ for $k > I_N$. Hence,

$$\sum_k \sqrt{\Pi(A_k)} \leq \sum_k \sqrt{\sum_{I_N \geq k} p_N}.$$

Consequently, if

$$\sum_k \sqrt{\bar{P}(M_k)} < \infty,$$

where $\bar{P}(M_k) = \sum_{N \geq M_k} p_N$ and $M_k = \min\{N : I_N \geq k\}$, then Hellinger consistency holds.

For example, if $I_N(\delta) = (c/\delta)^N$, for some $c > 0$ not depending on $\delta$, as is the case with Bernstein polynomials, then

$$M_k(\delta) = \lfloor \log k / \log(c/\delta) \rfloor$$

and, hence, we would wish that

$$\bar{P}(M_k(\delta)) < ak^{-2-r},$$

for some $r > 0$ and $a > 0$, for all large $k$ and for all $\delta > 0$. This holds if $\bar{P}(N) < a \exp(-N\psi)$ for all $\psi > 0$ for all large $N$, which holds if $N^{-1} \log \bar{P}(N) \to -\infty$ as $N \to \infty$.

6.4. *Random histogram.* Here we consider a random histogram model on $[0, 1]$ to illustrate Theorem 3. We take $m \in \{1, 2, \ldots\}$ with probability $\pi(m)$ and construct the random density function

$$f_m(x) = \sum_{k=1}^m w_{km} \mathbb{1}(a_{k-1m} < x < a_{km}),$$

where $w_{km} > 0$, $\sum_{k=1}^m w_{km} = m$ and $a_{km} = k/m$, $k = 0, 1, \ldots, m$. We will write $A_{km} = (a_{k-1m}, a_{km})$. We put $p_{km} = w_{km}/m$ and have a Dirichlet$(1 \ldots 1)$ prior for $p_m = (p_{1m}, \ldots, p_{mm})$. Then

$$f_n(x) = \sum_{m=1}^\infty f_{nm}(x)\pi(m|X^n),$$

where

$$f_{nm}(x) = \sum_{k=1}^{m} w_{kmn} \mathbb{1}(a_{k-1m} < x < a_{km})$$

and $w_{kmn} = \mathrm{E}(w_{km}|X^n)$. Also, for a nonnegative variable $Z$, $1/\mathrm{E}Z \leq \mathrm{E}1/Z$, so

$$\int f_0^2/f_n \leq \sum_m \pi(m|X^n) \int f_0^2/f_{nm}$$

and, therefore,

$$\mathrm{E}_{X^n} \int f_0^2/f_n \leq \sum_m \pi(m)\mathrm{E}_{X^n|m} \int f_0^2/f_{nm}.$$

Now $w_{kmn} = m(1+n_{km})/(m+n)$, where $n_{km} = \sum_{i=1}^{n} \mathbb{1}(X_i \in A_{km})$, and so

$$\mathrm{E}_{X^n|m} \int f_0^2/f_{nm} = \sum_{k=1}^{m} \int_{A_{km}} \mathrm{E}[(m+n)/\{m(1+n_{km})\}]f_0^2.$$

It is easy to show that

$$\mathrm{E}\{1/(1+n_{km})\} \leq 1/\{(1+n)F_0(A_{km})\}$$

with $n_{1m}\dots n_{mm} \sim \mathrm{mult}(n; F_0(A_{1m})\dots F_0(A_{mm}))$, and so

$$\mathrm{E}_{X^n} \int f_0^2/f_n \leq \sum_m \pi(m)\frac{m+n}{m(1+n)} \sum_{k=1}^{m} \int_{A_{km}} f_0^2/F_0(A_{km})$$

$$\leq \lambda \sum_m \frac{m+n}{1+n}\pi(m),$$

where $\lambda = \sup_x f_0(x)$, which we will assume to be finite. Therefore,

$$\sup_n \left\{ \mathrm{E}_{X^n} \int f_0^2/f_n \right\} < \infty$$

when $\sum_m m\pi(m) < \infty$ and so if the prior puts positive mass on all Kullback–Leibler neighborhoods of $f_0$, then $N^{-1}\sum_{n=1}^{N} D(f_{n-1}, f_0) \to 0$ a.s.

**7. Discussion.** As far as Hellinger consistency is concerned, the most fruitful sufficient conditions to date appear to be those involving the finiteness of sums of square roots of prior probabilities. Indeed, they improve on current sufficient conditions which have been published in the literature. In the case of Pólya trees the improvements are quite dramatic.

A framework for Kullback–Leibler consistency, which fits into a general framework including weak and Hellinger consistency, has been developed using martingales. Theorems 2, 2* and 3 suggest that the condition

$$\sum_n n^{-2} \mathrm{Var}\{\log(L_n/L_{n-1})\} < \infty$$

is highly significant and conditions under which this holds need to be understood.

Future work will investigate rates of convergence using the sums of square roots of prior probabilities approach. The basis for this is consideration of $\Pi^n(A_n)$, where $A_n = \{f : h(f, f_0) > \varepsilon_n\}$ and $\varepsilon_n \downarrow 0$. Following the proof of Theorem 4, we have

$$\Pi^n(A_n) \leq \exp(-n d_n)/\sqrt{I_n} \qquad \text{a.s.}$$

for all large $n$, for any sequence $d_n$ satisfying

$$\sum_n \exp\{-n(\gamma_n - d_n)\} K_n < \infty,$$

where

$$K_n = \sum_j \sqrt{\Pi(A_{nj})}$$

and $\{A_{nj}\}$ covers $A_n$ with $\delta_n$ size Hellinger balls and $\gamma_n = \varepsilon_n - \delta_n$. Putting $\delta_n = \varepsilon_n/2$ so $\gamma_n = \varepsilon_n/2$ seems appropriate here. Then, for example, lower bounds for $I_n$ in an a.s. sense are available from Shen and Wasserman ([2001](#)), using the $\rho_\alpha(f, f_0) = \alpha^{-1} \int \{(f_0/f)^\alpha - 1\} f_0$ metric, for $0 < \alpha \leq 1$. To find rates, it is required to understand $K_n$ which will be prior specific and involve a refinement of the work appearing in Section 6.

**Acknowledgments.** The author is grateful for the comments of a referee and Associate Editor whose suggestions resulted in a fundamental restructuring of the paper.

Department of Mathematical Sciences
University of Bath
Bath BA2 7AY
United Kingdom
E-mail: S.G.Walker@bath.ac.uk